\documentclass{ps}
\usepackage{amsmath,amsfonts,amssymb}
\usepackage{xcolor}
\usepackage[extra,safe]{tipa}
\usepackage[applemac]{inputenc}

\def \T{\mathbb{T}}
\def\my_c{c_\infty}

\def \z{{\mathbf{z}}}
 \def\bxi{{\boldsymbol{\xi}}}
 
 \def\bphi{{\boldsymbol{\phi}}}

 \def\btheta{{\boldsymbol{\theta}}}
\def\bTheta{{\boldsymbol{\Xi}}}

\def \cE{{\cal{E}}}

\def\lin{{\rm{linear}}}

\def\leftB{[\![}
\def\rightB{]\!]}
\newcommand \A[1]{{\bf (#1)}}

\def\R{{\mathbb{R}} }

\def\E{{\mathbb{E}}  }

\def\P{{\mathbb{P}}  }

\def\det{{\rm{det}}}

\def\bint#1^#2{\displaystyle{\int_{#1}^{#2}}}
\def\bsum#1^#2{\displaystyle{\sum_{#1}^{#2}}}

\def\xdt_#1{X_#1(\Delta t)}

\newtheorem{THM}{Theorem}[section]
\newtheorem{DEF}{Definition}[section]
\newtheorem{PROP}{Proposition}[section]

\newtheorem{COROL}{Corollary}[section]

\newtheorem{LEMME}{Lemma}[section]
\newtheorem{REM}{Remark}[section]

\newcommand{\mysection}{\setcounter{equation}{0} \section}

\def\diag{{{\rm diag}}}
\def\Tr{{{\rm Tr}}}

\def\0{{\mathbf{0}}}

\def\Cov{{\rm{Cov}}}

\def\X{{\mathbf{X}}}
\def\bx{{\mathbf{x}}}
\def\x{{\mathbf{x}}}
\def\by{{\mathbf{y}}}
\def\y{{\mathbf{y}}}

\def\gF{{\mathbf{F}}}
\def\bI{{\mathbf{I}}}

\def\G{{\mathbf{G}}}
\def\K{{\mathbf{K}}}
\def\gR{{\mathbf{R}}}

\def\gF{{\mathbf{F}}}
\def\gL{{\mathbf{L}}}

\begin{document}
%%-----------------------------
%%      the top matter
%%-----------------------------
\title{%Calder\'on-Zygmund estimates and associated 
Parametrix techniques and martingale problem for some degenerate Kolmogorov's equations} 
\author{S. Menozzi}\address{LPMA, Universit\'e Denis Diderot, 175 Rue du Chevaleret 75013 Paris, menozzi@math.jussieu.fr}

\date{\today}
\begin{abstract} We prove the uniqueness of the martingale problem associated to some degenerate %(weakly hypoelliptic) 
operators.
The key point is to exploit the strong parallel between the new technique introduced by Bass and Perkins  \cite{bass:perk:09} to prove uniqueness of the martingale problem in the framework of non degenerated elliptic operators and the Mc Kean and Singer \cite{mcke:sing:67} parametrix approach to the density expansion that has previously been extended to the degenerate  setting that we consider (see Delarue and Menozzi  \cite{dela:meno:10}). 
%controls deriving from the first step in the parametrix method of 
\end{abstract}
\maketitle

\mysection{Introduction}
\label{INTRO}
\subsection{Martingale problem and parametrix techniques}
The martingale approach turns out to be particularly useful when trying to get uniqueness results for the stochastic process corresponding to an operator. 
In a recent work, R. Bass and E. Perkins \cite{bass:perk:09} introduced in the framework of non-degenerated, non-divergence, time-homogeneous operators a new technique to prove uniqueness of the associated martingale problem. Precisely, for an operator of the form
\begin{equation}
\label{EQ1}
Lf(x)=\frac{1}2\Tr (a(x) D_x^2f(x)), \ f\in C_0^2(\R^d,\R),\ x\in \R^d,
\end{equation}
the authors prove uniqueness provided $a$ is uniformly elliptic, bounded and uniformly $\eta $-H\"older continuous in space ($\eta \in (0,1]$), i.e. there exists $C>0 $ s.t. for all $(x,y)\in \R^d,\ |a(x)-a(y)|\le C(1\wedge |x-y|^\eta) $. That is, for a given starting point $x\in \R^d$, there exists a unique probability measure $\P$ on $C(\R^+ ,\R^d) $ s.t. denoting by $(X_t)_{t\ge 0}$ the canonical process, $\P[X_0=x]=1 $ and for every $f\in C_0^2(\R^d,\R)$, $f(X_t)-f(x)-\int_0^tLf(X_s) ds $ is a $\P$-martingale.

In the indicated framework, this result can be derived from the more involved Calder\'on-Zygmund like $L^p$ estimates  established by Stroock and Varadhan \cite{stro:vara:79}, that only require continuity of the diffusion matrix $a$, or from a more analytical viewpoint from some appropriate Schauder estimates, see e.g. Friedman \cite{frie:64}.

Anyhow, the technique introduced in \cite{bass:perk:09} can be related with the first step of Gaussian approximation of the parametrix expansion of the fundamental solution of \eqref{EQ1} developed by McKean and Singer \cite{mcke:sing:67} that we now shortly describe. Suppose first that, additionally to the previous assumptions of ellipticity, boundedness and uniform H\"older continuity, the diffusion coefficient $a$ is smooth (say $C^\infty(\R^d,\R)$). Thus, the fundamental solution $p(s,t,x,y)$ of \eqref{EQ1} exists and is smooth for $t>s$, see e.g \cite{frie:64}. Precisely, we have:
\begin{eqnarray*}
\partial_t p(s,t,x,y)=L^* p(s,t,x,y),\ t>s,\ (x,y)\in (\R^d)^2,\ p(s,t,x,.)\underset{t\downarrow s}{\longrightarrow} \delta_{x}(.),
\end{eqnarray*}
where $L^*$ stands for the adjoint of $L$ and acts on the $y$ variable.
For fixed starting and final points $x,y\in \R^d$ and a given final time $t>0$, in order to estimate $p(0,t,x,y) $, one introduces the Gaussian process $\tilde X_u^y=x+\sigma(y)W_{u-s}, \ u\in [s,t],s\le t $, where $(W_u)_{u\in [0,t-s]} $ is a standard $d $-dimensional Brownian motion and $\sigma\sigma^*(y)=a(y) $. Observe that the coefficient of $\tilde X^y$ is frozen here at the point where we consider the density. Denote by $\tilde p^y(s,t,x,.)$ the density of  $\tilde X^y $ at time $t$ starting from $x$ at time $s$, and for $\varphi \in C_0^2(\R^d,\R)$, define by $\tilde L^y \varphi(x)=\frac12 \Tr(a(y)D_{x}^2\varphi(x) )$ its generator. The density $\tilde p^y(s,t,x,.)$ satisfies the Kolmogorov equation:
\begin{eqnarray*}
\partial_s \tilde p^y(s,t,x,z)=-\tilde L^y \tilde p^y(s,t,x,z),\ s<t,\ (x,z)\in (\R^d)^2,\ \tilde p^y(s,t,.,z)\underset{s\uparrow t}{\longrightarrow } \delta_z(.),
\end{eqnarray*}
where $\tilde L^y $ acts here on the $x$ variable.
Take now $z=y$ in the above equation. By formal derivation and the previous Kolmogorov equations we obtain:
\begin{eqnarray}
p(0,t,x,y)-\tilde p^y(0,t,x,y)&=&\int_0^t ds \partial_s \int_{\R^d} p(0,s,x,w)\tilde p^y(s,t,w,y)dw\nonumber\\
                                      &=&\int_{0}^t ds\bint{\R^d}^{}\left( L^* p(0,s,x,w)\tilde p^y(s,t,w,y)-p(0,s,x,w)\tilde L^y\tilde p^y(s,t,w,y)\right) dw \nonumber\\
                                      &=&\int_0^t ds\bint{\R^d}^{} p(0,s,x,w) (L-\tilde L^y) \tilde p^y(s,t,w,y) dw\nonumber \\
                                      &=&\int_0^t ds\bint{\R^d}^{} p(0,s,x,w) H  (s,t,w,y) dw:=p\otimes  H(0,t,x,y), \label{PARAM_P}
\end{eqnarray}
where %$L^*$ denotes the formal adjoint of $L$ and $\tilde L $ the generator of $\tilde X$ and 
$\otimes $ denotes a time-space convolution. Observe that $H(s,t,w,y):=(a(w)-a(y))D_z^2 \tilde p^y(s,t,w,y) $. From direct computations, there exist $c,C>0$ (depending on $d$, the uniform ellipticity constant and the $L_\infty$ bound of $a$) s.t. $|D_z^2\tilde p^y(s,t,w,y)| \le \frac{C}{(t-s)^{d/2+1}}\exp(-c\frac{|w-y|^2}{t-s})$.
%Hence, there exists $\tilde c,\ $ $\le \frac{\tilde C}{(t-s)^{d/2+(1-\eta/2)}}\exp(-\tilde c\frac{|w-y|^2}{t-s})$. 
The previous uniform H\"older continuity assumption on $a$ is therefore a sufficient (and quite sharp) condition  to remove the time-singularity in $H$. The idea of the parametrix expansion is then to proceed in \eqref{PARAM_P} by applying the same freezing technique to $p(
0,s,x,w)$ introducing the density $\tilde p^w(0,s,x,.) $ of the process with coefficients frozen at point $z$. One eventually gets the formal expansion 
\begin{eqnarray}
\label{formal_Series}
p(0,t,x,y)=  \tilde p(0,t,x,y)+\bsum{k\ge 1}^{}\tilde p\otimes H^{\otimes k}(0,t,x,y),
\end{eqnarray}
where $H^{\otimes k},\ k\ge 1 $, stands for the iterated convolutions of $H$, and $\forall (s,t,z,y)\in (\R^+)^2\times \R^d,\ \tilde p(s,t,z,y):=\tilde p^y(s,t,z,y) $. The H\"older continuity gives that $H$ is a ``smoothing" kernel in the sense that there exist $c,C>0$ (with the same previous dependence) s.t.  $\ |\tilde p\otimes H^{\otimes k}(s,t,z,y)|\le C^{k+1}(t-s)^{k\eta/2} \prod_{i=1}^{k+1}B\left(1+\frac{(i-1)\eta}{2},\frac \eta 2\right) (t-s)^{-d/2 }\exp(-c\frac{|y-z|^2}{t-s})$, where $B(m,n)=\int_{0}^1 s^{m-1}(1-s)^{n-1} d s $ stands for the $\beta $ function. From this estimate, equation \eqref{formal_Series} and the asymptotics of the $\beta$ function, one directly gets the Gaussian upper bound over a compact time interval. Namely for all $T>0$, there exist constants $c,C>0$ s.t.
\begin{equation}
\label{ARONSON_UE}
\forall 0\le s<t\le T,\ \forall (x,y)\in (\R^d)^2,\ p(s,t,x,y)\le \frac{C}{(t-s)^{d/2}} \exp(-c\frac{|y-x|^2}{t-s}),
\end{equation}
with $c,C$ depending on $d$, the uniform ellipticity constant and $L_\infty$ bound of $a$ and $C$ depending on $T$ as well. 
We refer to Konakov and Mammen \cite{kona:mamm:00} for details in this framework.

Up to now we supposed $a$ was smooth in order to guarantee the existence of the density and justify the formal derivation in \eqref{PARAM_P}. On the other hand, the r.h.s. of \eqref{formal_Series} can be defined without additional smoothness on $a$ than uniform $\eta $-H\"older continuity. The Gaussian upper bound \eqref{ARONSON_UE} also only depends on the H\"older regularity of $a$. A natural question is to know whether the r.h.s. of \eqref{formal_Series} corresponds to the density of some stochastic differential equation under the sole assumptions of uniform ellipticity, boundedness and H\"older continuity on $a$. A positive answer is given by the uniqueness of the martingale problem associated to \eqref{EQ1}. Indeed, considering a sequence of equations with mollified coefficients, we derive from convergence in law, the Radon-Nikodym theorem and \eqref{ARONSON_UE}  that the unique weak solution of $dX_t=\sigma(X_t)dW_t $ associated to $L$ admits a density that satisfies the previous Gaussian bound. It is actually remarkable that the uniqueness of the martingale problem can be proved using exactly the smoothing properties of the previous kernel $H$. That is what was achieved by Bass and Perkins \cite{bass:perk:09} in the framework we described and it is the main purpose of this note in a degenerate setting.

To conclude this paragraph, let us emphasize that the previous parametrix approach has been used in various contexts. It turns out to be particularly well suited to the approximation of the underlying processes by Markov chains, see Konakov and Mammen \cite{kona:mamm:00,kona:mamm:02} for the non degenerate continuous case or \cite{kona:meno:10} for the approximation of stable driven SDEs. On the other hand, recently, we used this technique to give a local limit theorem for the Markov chain approximation of a Langevin process \cite{kona:meno:molc:10} or two-sided bounds of some more general degenerated hypoelliptic operators \cite{dela:meno:10}.
In particular, in both works, we have an unbounded drift term.
The unboundedness of the first order term imposes a more subtle strategy than the previous one for the choice of the frozen Gaussian density. Namely, one has to take into consideration in the frozen process the ``geometry" of the deterministic differential equation 
associated to the first order terms of the operator. This will be thoroughly explained in the next section. Anyhow, the strategy of the previous articles allows to extend the technique of Bass and Perkins to prove uniqueness of the martingale problem for some degenerate operators with unbounded coefficients. %This is the main purpose of the current work.

\subsection{Statement of the Problem and Main Results}
 
Consider the following system of Stochastic Differential Equations (SDEs in short)
\begin{equation}
\label{SYST}
\begin{array}{l}
\displaystyle dX_t^1 = F_1(t,X_t^1,\dots,X_t^n) dt + \sigma(t,X_t^1,\dots,X_t^n) dW_t,
\\
\displaystyle dX_t^2 = F_2(t,X_t^1,\dots,X_t^n) dt,
\\
\displaystyle dX_t^3 = F_3(t,X_t^2,\dots,X_t^n) dt,
\\
\displaystyle \cdots
\\
\displaystyle dX_t^n = F_n(t,X_t^{n-1},X_t^n) dt,
\end{array}
\quad t \geq 0,
\end{equation}
$(W_t)_{t \geq 0}$ standing for a $d$-dimensional Brownian motion, and each $(X_t^i)_{t \geq 0}$, $1 \leq i \leq n$, being $\R^d$-valued as well.
%%%%% Il faut aussi dire un mot des applications

From the applicative viewpoint, systems of type \eqref{SYST} appear in many fields. Let us for instance mention for $n=2$ stochastic Hamiltonian systems (see e.g. Soize \cite{soiz:94} for a general overview or Talay \cite{tala:02} and H\'erau and Nier \cite{hera:nier:04} for convergence to equilibrium). 
Again for $n=2$,  the above dynamics is used in mathematical finance to price Asian options (see for example 
\cite{baru:poli:vesp:01}). For $n \geq 2$, it appears in heat conduction models (see e.g. Eckmann et al. 
\cite{eckm:99} and Rey-Bellet and Thomas \cite{reyb:thom:00} when the chain is forced by two heat baths).

In what follows, we denote a quantity in $\R^{nd}$ by a bold letter: i.e. $\0$, stands for zero in $\R^{nd}$ and the solution $(X_t^1,\dots,X_t^n)_{t \geq 0}$ 
to \eqref{SYST} is denoted by $({\mathbf X}_t)_{t \geq 0}$. Introducing the embedding matrix $B$ from $\R^d$ into $\R^{nd} $, i.e. $B= (I_d , 0, \dots, 0)^*$, where ``$*$'' stands for the transpose, we rewrite \eqref{SYST} in the shortened form 
\begin{equation*}
d{\mathbf X}_t = {\mathbf F}(t,{\mathbf X}_t) + B \sigma(t,{\mathbf X}_t) dW_t,
\end{equation*}
where ${\mathbf F}=(F_1,\dots,F_n)$ is an $\R^{nd}$-valued function.
Moreover, for $\bx=(x_1,\dots,x_n) \in (\R^d)^n$, we set $\bx^{i,n} = (x_i,\dots,x_n)
\in (\R^d)^{n-i+1}$.

We introduce the following assumptions:
\begin{trivlist}
\item[\A{R-$\eta$}]  The functions $(F_ i)_{i\in \leftB 1,n\rightB}$ are uniformly Lipschitz continuous with constant $\kappa>0 $ (alternatively we can suppose for $i=1$ that the drift of the non degenerated component $F_1 $ is measurable and bounded by $\kappa$). The diffusion matrix $(a(t,.))_{t\ge 0} $ is uniformly $\eta $-H\"older continuous in space with constant $\kappa $, i.e.
$$  \forall t\ge 0,\ \sup_{(\x,\y) \in \R^{nd},\ \x\neq \y}\frac{|a(t,\x)-a(t,\y)|}{|\x-\y|^\eta} \le \kappa.$$
%for all $T>0$ and $x\in \R^{nd}$.
\item[\A{UE}] There exists  $\Lambda \ge 1,\ \forall t\ge 0, \x\in \R^{nd},\ \xi \in \R^d, \ \Lambda^{-1}|\xi|^2\le \langle a(t,\x) \xi,\xi\rangle \le \Lambda |\xi|^2 $. 
\item[\A{ND-$\eta $}] For each integer $2 \leq i \leq n$, $(t,(x_i,\dots,x_n)) \in \R_+ \times \R^{(n-i+1)d}$, the function $x_{i-1} \in \R^{d} \mapsto F_i(t,\bx^{i-1 , n})$
is continuously differentiable, the derivative,
%denoted by 
$(t,\x^{i-1,n}) \in \R_+ \times \R^{(n-i+2)d}\mapsto
D_{x_{i-1}} F_i(t,\x^{i-1,n})$, being $\eta $-H\"older continuous %in the first space variable $x_{i-1} $ % on zappe pour l'instant cette hypothse
with constant $\kappa $.
There exists a closed convex subset ${\mathcal E}_{i-1} \subset GL_d(\R)$
(set of invertible $d \times d$  matrices)
 s.t., for all
$t \geq 0$ and $(x_{i-1},\dots,x_n) \in \R^{(n-i+2)d}$, the matrix $D_{x_{i-1}} F_i(t,\x^{i-1 ,n})$
belongs to ${\mathcal E}_{i-1}$.
For example, ${\mathcal E}_i$, $1 \leq i \leq n-1$, may be a closed ball
included in $GL_d(\R)$, which is an open set.\\
 \end{trivlist}

Assumptions \A{UE}, \A{ND-$\eta$} can be seen as a kind of (weak) H\"ormander condition. They allow to transmit the non degenerate noise of the first component to the other ones. Also, the particular structure of $F(t,.)=(F_1(t,.),\cdots, F_n(t,.)$ yields that the $i^{{\rm th}} $ component has intrinsic time scale $(2i-1)/2, i\in \leftB 1, n\rightB $. 
We notice that the coefficients may be irregular in time.  The last part of Assumption \A{ND-$\eta$} will be explained in Section \ref{SEC_GAUSS}. We say that assumption \A{A-$\eta$} is satisfied if \A{R-$\eta$}, \A{UE}, \A{ND-$\eta$} hold.

Under \A{A-$\eta$}, we established in \cite{dela:meno:10} Gaussian Aronson like estimates for the density of \eqref{SYST} over compact time interval $[0,T]$, for $\eta>1/2$. Precisely, we proved that
the unique weak solution of \eqref{SYST} admits a density that satisfies that for all $T>0, \exists C:=C(T,$\A{A-$\eta$}$)$ s.t. $\forall (t,\x,\y)\in (0,T]\times \R^{nd}\times \R^{nd}$:
\begin{equation}
\label{BD_ARONSON_DEG}
C^{-1}t^{-n^2d/2}\exp\left(-Ct|\T_t^{-1}(\btheta_t(\x)-\y|^2 \right) \le p(t,\x,\y)\le Ct^{-n^2d/2}\exp\left(-C^{-1}t|\T_t^{-1}(\btheta_t(\x)-\y|^2 \right),
\end{equation}
where $\T_t:=\diag( (t^i I_d)_{i\in \leftB 1,n \rightB })$ is a scale matrix and $\overset{.}{\btheta}_t(x)= \gF(t,\btheta_t(\x)), \btheta_0(\x)=\x $.

To derive \eqref{BD_ARONSON_DEG}, we proceeded using a ``formal" parametrix expansion considering a sequence of equations with smooth coefficients for which H\"ormander's theorem guaranteed the existence of the density, see. e.g. H\"ormander \cite{horm:67} or Norris \cite{norr:86}.  Anyhow, as in the previous paragraph, our estimates did not depend on the derivatives of the mollified coefficients but only on the $\eta $-H\"older continuity assumed in \A{A-$\eta $}. Anyhow, to pass to the limit following the previously described procedure, some uniqueness in law is needed. Using the comparison principle for viscosity solutions of fully non-linear PDEs,  see Ishii and Lions \cite{ishi:lion:90}, we managed to obtain the bounds under \A{A-$\eta$},  $\eta >1/2$.  However, the viscosity approach totally ignores the smoothing effects of the heat kernel and is not a ``natural technique" to derive uniqueness in law.

Introduce the generator of \eqref{SYST}:
\begin{equation}
\label{THE_GEN}
\forall \varphi \in C_0^2(\R^+\times \R^{nd},\R), \ \forall (t,\x) \in \R^+\times \R^{nd},\ L_t\varphi (t,\x)=\langle \gF(t,\x), {\mathbf D}_\x \varphi(t,\x)\rangle+\frac 12 \Tr( a(t,\x) D_{\x_1}^2 \varphi(t,\x)) . 
\end{equation}

Adapting the technique of Bass and Perkins \cite{bass:perk:09} we obtain the following results.
\begin{THM}
\label{THM_M}
For $\eta \in (0,1] $, under \A{A-$\eta $} the martingale problem associated to $L$ in \eqref{THE_GEN} is well-posed. In particular, weak uniqueness in law holds for the SDE \eqref{SYST}.
\end{THM}
As a bypass product we derive from \cite{dela:meno:10}  the following:
\begin{COROL}
\label{COROL_THM1} For $\eta $ in $(0,1] $, under \A{A-$\eta $},  the unique weak solution of \eqref{SYST} admits for all $t>0 $ a density that satisfies the Aronson like bounds of equation \eqref{BD_ARONSON_DEG}.
\end{COROL}

\mysection{Choice of the reference Gaussian process for the parametrix}
In this section we describe the Gaussian processes that will be involved in the study of the martingale problem and that have been previously involved in the parametrix expansions of \cite{dela:meno:10}. We first introduce in Section \ref{SEC_GAUSS} a class of degenerate linear stochastic differential equations that admit a density satisfying bounds similar to those of equation \eqref{BD_ARONSON_DEG}. We then specify, how to properly linearize the dynamics of \eqref{SYST} so that the linearized equations belong to the class considered in Section \ref{SEC_GAUSS}.

\subsection{Some estimates on degenerate Gaussian processes with linear drift}
\label{SEC_GAUSS}

Introduce the stochastic differential equation:
\begin{equation}
\label{EQ_LIN}
d\G_t=\gL_t \G_t +B \Sigma_t dW_t
\end{equation}
where ${\mathbf L}_t \bx
= (0,\alpha_t^1x_1,\dots,\alpha_t^{n-1} x_{n-1})^* + {\mathbf U}_t \bx$, and 
${\mathbf U}_t\in\R^{nd} \otimes \R^{nd}$ is an  ``upper triangular'' block matrix with zero entries on its first $d$ rows. 
We suppose that the coefficients satisfy the following assumption \A{A$^{{\rm linear}}$}:
\begin{trivlist}
\item[-] The diffusion coefficient $A_t:=\Sigma_t\Sigma_t^*, \ t\ge 0$, is uniformly elliptic and bounded, i.e. $\exists \Lambda\ge 1, \ \forall \xi \in \R^d,\ \Lambda^{-1}|\xi|^2\le \langle A_t \xi , \xi\rangle \le \Lambda |\xi|^2  $. 
\item[-] 
For each $i \in \leftB 1,n-1\rightB $, there exists a closed convex subset ${\mathcal E}_{i} \subset GL_d(\R)$ s.t. for all
$t \geq 0$, the matrix $\alpha_t^{i}$ belongs to ${\mathcal E}_{i}$.
\end{trivlist}

Denoting by $(\gR(s,t))_{0\le s,t} $ the resolvent associated to $(\gL_t)_{t\ge 0} $, i.e. $\partial_t \gR(t,s)=\gL_t \gR(t,s),\ \gR(s,s)=\bI_{nd}$, we have 
for $0\le s<t, \x\in \R^{nd} $,  $\G_t^{s,x}=\gR(t,s)x+\int_s^t  \gR(t,u) B\Sigma_u dW_u $ and $\Cov[\G_t^{s,x}]=\int_s^t \gR(t,u) B A_u B^*\gR(t,u)^* du :=\K(s,t)$.

From Propositions 3.1 and 3.4 in \cite{dela:meno:10}, the family $(\K(s,t))_{t\in (s,T]} $ of covariance matrices associated to the Gaussian process $(\G_t^{s,\x})_{t\in (s,T]} $ satisfies, under \A{A$^{{\rm linear}}$}, a ``good scaling property" in the following sense:
\begin{DEF}[Good scaling property]
\label{GSP}
Fix $T>0$. We say that a family $(\K(s,t))_{t\in [s,T]}, s\in [0,T) $ of $\R^{nd}\otimes \R^{nd} $ matrices satisfies a good scaling property with constant $C\ge 1$ (see also Definition 3.2 and Proposition 3.4 of \cite{dela:meno:10}) if for all $(t,s)\in (\R^+)^2,\ 0<  t-s\le T, \ \forall \y \in \R^{nd},\ C^{-1}(t-s)^{-1}|\T_{t-s} \y|^2\le \langle \K(s,t)\y,\y\rangle\le C  (t-s)^{-1}|\T_{t-s}\y|^2 $. 
\end{DEF}
Precisely the family $(\K(s,t))_{t\in (s,T]} $ satisfies under \A{A$^\lin$} a good scaling property with constant $C:=C(T,$\A{A$^\lin $}$)$.
\begin{REM}
We point out that it is precisely the second assumption of \A{A$^\lin$} concerning the existence of convex subsets $(\cE_i)_{i\in \leftB 1,n-1\rightB} $ of $GL_d(\R) $ that guarantees the good scaling property (see Propositions 3.1 and 3.4 in \cite{dela:meno:10} for details).
\end{REM}

The density at time $t>s$ in $\y\in \R^{nd} $ of  $\G_t^{s,\x} $ writes
\begin{equation}
\label{DENS_LIN}
q(s,t,\x,\y)=\frac{1}{(2\pi)^{nd/2}\det (\K_{s,t})^{1/2}}\exp(-\frac 12\langle \K(s,t)^{-1}(\gR(t,s)\x-\y), \gR(t,s)\x-\y \rangle).%\exp\left(-\frac 12 \langle \K(s,t)^{-1}(\gR(t,s)\x-\y), \gR(t,s)\x-\y \rangle \right)
\end{equation}
%where $I(s,t,\x,\y):=\inf\{\int_s^t|u_v|^2dv,\ u\in L^2([s,t],\R^d), \bvarphi(s)=\x, \ \bvarphi(t)=\y,\ \overset{.}{\bvarphi}(v)=\gL_v \bvarphi(v)+Bu_v, \ v\in [s,t) \} =\langle \K(s,t)^{-1}(\gR(t,s)\x-\y), \gR(t,s)\x-\y \rangle$ and $\K(s,t)=\int_s^t \gR(t,u) B A_u B^*\gR(t,u)^* du =\Cov[\G_t^{s,x}]$ which is also the Gram matrix associated with $I$. We refer to Theorem 1.11, Chapter 1 of Coron \cite{coro:07} for a justification a the last equality, which is the explicit solution of a linear controllability problem. 

Since under \A{A$^\lin$}, $(\K(s,t))_{t\in (s,T]} $ satisfies a good scaling property in the sense of Definition \ref{GSP}, we then derive from  \eqref{DENS_LIN}:
\begin{PROP}
\label{CT_AR_L}
Under \A{A$^{\lin}$}, for all $T>0$ there exists a constant $C_{\ref{CT_AR_L}}:=C_{\ref{CT_AR_L}}(T,\A{A^{\lin}})\ge 1$ s.t. :
\begin{eqnarray}
&&C_{\ref{CT_AR_L}}^{-1} (t-s)^{n^2d/2}\exp(-C_{\ref{CT_AR_L}}(t-s)|\T_{t-s}^{-1}(\gR(t,s)\x-\y)|^2) \le     q(s,t,\x,\y)\nonumber\\
  &&\hspace*{3cm}\le C_{\ref{CT_AR_L}} (t-s)^{n^2d/2}\exp(-C_{\ref{CT_AR_L}}^{-1} (t-s)|\T_{t-s}^{-1}(\gR(t,s)\x-\y)|^2). \label{EQ_LIN_A}
\end{eqnarray}
\end{PROP}
 This means that the off-diagonal bound of Gaussian processes with dynamics \eqref{EQ_LIN} and fulfilling \A{A$^{\lin}$} is homogeneous to the square of the difference between the final point $\y $ and $\gR(t,s)\x $ (which corresponds to the transport of the initial condition by the deterministic system deriving from \eqref{EQ_LIN}, that is $\overset {.}{\btheta}_t=\gL_t \btheta_t $) rescaled by the intrinsic time-scale of each component. We here recall that the component $i\in \leftB 1,d \rightB$  has characteristic time scale $t^{(2i-1)/2}$.

\subsection{Linearization of the initial dynamics and associated estimates}
\label{SEC_FRO}
The crucial feature of the parametrix method described in the introduction was to choose a ``good" process to approximate the density of the diffusion.
In the uniformly elliptic case, with bounded coefficients, one could take, as a first approximation, the Gaussian process with coefficients frozen in space at the fixed final spatial point where we wanted to estimate the density. The choice is natural since it makes the kernel $H$ (defined in \eqref{PARAM_P}) ``compatible" with the bounds of the frozen density. It is precisely the off diagonal term in $\exp(-c|x-y|^2/(t-s)) $ that allows to equilibrate the singularity in $|x-y|^\eta/(t-s) $ coming from the second order spatial derivatives.  In their work, \cite{bass:perk:09}, Bass and Perkins exactly exploited the specific behavior of the singular kernel $H$  which has an integrable singularity in time at 0 (see their Proposition 2.3), to derive uniqueness of the martingale problem in the non-degenerate time-homogeneous framework. This approach provides a natural link between parametrix expansions and the study of martingale problems for uniformly H\"older continuous coefficients.

Parametrix expansions, to derive density estimates on systems of the form \eqref{SYST}, have been discussed in \cite{dela:meno:10}. We thus have in the current degenerate framework of assumption \A{A-$\eta $} some natural candidate defined below. The key idea is to consider a ``degenerate" Gaussian process whose density anyhow has a specific ``off-diagonal" behavior similar to the one exhibited in equation \eqref{EQ_LIN_A} and to choose the freezing process in order that the singularity deriving  from $H$ is still compatible with the  ``off-diagonal" bound in the sense that it will be sufficient to remove the time-singularity.
%Intuitively the ``off-diagonal" term corresponds to the weight associated to a deviation due to randomness from the behavior of the deterministic equation.  Control theory provides a precise formulation of this intuitive fact.  $\exp(-cI(T,\x,\y) )$ where $I(T,\x,\y):=\inf\{\int_0^T|u_s|^2, \bvarphi(0)=\x, \ \bvarphi(T)=\y,\ \overset{.}{\bvarphi}(t)=\gF(t,\bvarphi(t))+Bu_t \} $ is the energy to make the deteriministic system go from $\x$ at time $0$ to $\y$ at time $T$. This fact was thoroughly used in \cite{dela:meno:10}.
%%%%% A developper

 We follow the same line of reasoning in our current framework.

%In the current degenerate case, since we consider as well unbounded coefficients, the crucial question concerns the choice of the Gaussian process
%that will serve as first approximation of the dynamics.  
For fixed parameters $T>0, \y \in \R^{nd}$, introduce the linear equation:
\begin{equation}
\label{eq:F:linear:eq:upper}
\begin{split}
d \tilde{\mathbf X}_t^{T,\y} &= \bigl[ {\mathbf F}(t,{ \btheta}_{t,T}(\y)) 
+ D{\mathbf F}(t,{ \btheta}_{t,T}(\y))\bigl(\tilde{\mathbf X}_t^{T,\y} - {\btheta}_{t,T}(\y) \bigr) \bigr] dt
\\
&\hspace{5pt} +  B \sigma(t,{\btheta}_{t,T}(\y)) dW_t, \quad 0 \leq t \leq T,
\end{split}
\end{equation}
where $({\btheta}_{t,T}(\y))_{t  \geq 0}$ solves the ODE
$[d/dt]{\btheta}_{t,T}(\y) = {\mathbf F}(t,{\btheta}_{t,T}(\by))$, $t \geq 0$,
with the boundary condition ${\btheta}_{T,T}(\y)=\y$ and $\forall (t,\x)\in [0,T]\times \R^{nd},\ D\gF(t,\x)=\left (\begin{array}{ccccc}0 & \cdots & \cdots &\cdots  & 0\\
D_{\x_1}\gF_2(t,\x) & 0 &\cdots &\cdots &0\\
0 & D_{\x_2} \gF_3(t,\x)& 0& 0 &\vdots\\
\vdots &  0                      & \ddots & \vdots\\
0 &\cdots &     0      & D_{\x_{n-1}}\gF_n(t,\x) & 0
\end{array}\right) 
 $ is the subdiagonal of the Jacobian matrix ${\mathbf D_\x \gF} $. 
%To simplify, we do not specify the dependence of $\tilde{\mathbf X}$ on $(T,\by)$ and 
Write $\tilde p^{T,\y}(t,T,\x,.) $ for the density of $\tilde \X_T^{T,\y} $ starting from $\x$ at time $t$.

The deterministic ODE associated with 
$\tilde{\mathbf X}^{T,\y}$ has the form
\begin{equation}
\label{eq:F:240409:3}
\frac{d}{dt} \tilde{\bphi}_t = {\mathbf F}(t,\btheta_{t,T}(\by))
+ D {\mathbf F}(t,\btheta_{t,T}(\by))[\tilde{\bphi}_t
- \btheta_{t,T}(\by)], \quad t \geq 0.
\end{equation}
We denote by $(\tilde{\btheta}_{t,s}^{T,\by})_{s,t \geq 0}$
the associated flow, i.e. $\tilde{\btheta}_{t,s}^{T,\by}(\bx)$ is the value of $\tilde{\bphi}_t$ when $\tilde{\bphi}_s = \bx$. It is affine:
\begin{equation}
\label{AFFINE}
\begin{split}
\tilde{\btheta}_{t,s}^{T,\by}(\bx)&=\tilde \gR^{T,\y}(t,s)\x
 \\
 &\hspace{5pt}
+\int_{s}^t \tilde \gR^{T,\y}(t,u)\bigl({\mathbf F}(u,\btheta_{u,T}(\by))
- D {\mathbf F}(u,\btheta_{u,T}(\by)) \btheta_{u,T}(\by)\bigr) du.
\end{split}
\end{equation}
Above, $(\tilde \gR^{T,\y}(t,s))_{s,t \geq 0}$ stands for the resolvent associated with the matrices $(D {\mathbf F}(t,\btheta_{t,T}(\by)))_{t \geq 0} $.
%We then know
%\begin{equation}
%\label{eq:F:270409:1}
%-\ln \bigl[\tilde{p}(t,\varepsilon;\by,\bx')\bigr] + n^2d/2 
%\ln(\varepsilon-t)
%\approx
%\varepsilon |{\mathbb T}_{\varepsilon-t}^{-1}
%(\phi_{t,\varepsilon}(\by) - \bx')|^2 + O(1).
%\end{equation}
%(We give a precise meaning to the symbol $\approx$ in the sequel.)

We now claim
\begin{LEMME}
\label{eq:F:240409:2}
Let $T_0>0$ be fixed. There exists a constant $C_{\ref{eq:F:240409:2}}\geq 1$,
depending on ({\bf A}) and $T_0$ such that, for any $t \in [0,T),\ T\le T_0$ and 
$\bx,\by \in \R^{nd}$, 
\begin{equation*}
C_{\ref{eq:F:240409:2}}^{-1}  \bigl|{\mathbb T}_{T-t}^{-1}
\bigl[\bx - \btheta_{t,T}(\by) \bigr] \bigr| \leq 
 \bigl|{\mathbb T}_{T-t}^{-1}
\bigl[\tilde{\btheta}_{T,t}^{T,\by}(\bx) - \by \bigr] \bigr|
\leq C_{\ref{eq:F:240409:2}}  
\bigl|{\mathbb T}_{T-t}^{-1}
\bigl[\bx - \btheta_{t,T}(\by)\bigr]\bigr|. 
\end{equation*}
\end{LEMME}
This means that we can compare the rescaled ``forward" transport of the initial condition $\x$ from $t$ to $T$ by the linear flow and the rescaled ``backward" transport from $T$ to $t$ 
of the final point $\y$ by the original deterministic differential dynamics. We refer to Lemma 5.3 of \cite{dela:meno:10} for a proof.

Furthermore, under \A{A-$\eta $} we have that $D\gF(t,\btheta_{t,T}(\y))$ satisfies \A{A$^\lin$}. We thus derive from Lemma \ref{eq:F:240409:2} and a direct extension of Proposition \ref{CT_AR_L} (the mean of $\tilde{\mathbf X}_T^{T,\y} $ starting from $\x $ at time $t$ being $\tilde \btheta_{T,t}^{T,\y}(\x) $)  the following result.
\begin{LEMME}
\label{lem:F:040509:1}
Let $T_0>0$ be fixed. There exists 
a constant $C_{\ref{lem:F:040509:1}} >0$,
depending on ({\bf A}-$\eta$) and $T_0$
%only (and not on $T$), 
such that, for all $0 \leq t <T\le T_0$
and 
$\bx, \by \in \R^{nd}$, 
\begin{equation*}
\tilde{p}^{T,\y}(t,T,\bx,\by)
\leq C_{\ref{lem:F:040509:1}}
g_{C_{\ref{lem:F:040509:1}},T-t}\bigl(\bx - \btheta_{t,T}(\by)\bigr),
\end{equation*}
where for all $a>0,t>0,\ g_{a,t}(\by) = 
t^{-n^2\frac{d}{2}}
\exp(-a^{-1} t
|{\mathbb T}_{t}^{-1} \by|^2), \quad a,t>0, \ \by \in \R^{nd}$.
\end{LEMME}
%These Lemmas are proved  and 5.4 of that reference).
For all $ 0\le s<t\le T,\ \z,\y\in \R^{nd}$, define now the kernel $H$ as:
\begin{eqnarray}
\label{DEF_KERN_H}
 H(s,t,\z,\y)&=&\{ \langle \gF(s,\z) ,{\mathbf D}_\z  \tilde p^{t,\y}(s,t,\z,\y)\rangle+\frac 12 \Tr(a(s,\z)D_{\z_1}^2\tilde p^{t,\y}(s,t,\z,\y))\}-\nonumber\\
&&\{ \langle \gF(s,\btheta_{s,t}(\y))+D\gF(s,\btheta_{s,t}(\y))(\z-\btheta_{s,t}(\y)) ,{\mathbf D}_\z\tilde p^{t,\y}(s,t,\z,\y)\rangle \nonumber\\
&&+\frac 12 \Tr(a(s,\btheta_{s,t}(\y) )D_{\z_1}^2 \tilde p^{t,\y}(s,t,\z,\y) )\}  %\nonumber\\
 %&
 =
 %&
 (L_{s,\z}-\tilde L_{s,\z}^{t,\y})\tilde p^{t,\y}(s,t,\z,\y),
\end{eqnarray}
where $L$ is the generator of the initial diffusion \eqref{SYST} defined in \eqref{THE_GEN}, $\tilde L_{s,\z}^{t,\y} $ and $\tilde p^{t,\y}(s,t,\z,.)$ respectively stand for the generator at time $s$ and the density at time $t$ of $\tilde \X^{t,\y},\ \tilde \X_{s}^{t,\y}=\z $ with coefficients ``frozen" w.r.t. $t,\y$. The lower script in $\z$ is to emphasize that $\z $ is the differentiation parameter in the operator. 
%$$\forall \varphi\in C_0^2,\ \tilde L_{}^{s,\z} $$
%Observe in particular that $L_{s,\z}=\tilde L_{s,z}^{s,\z} $.

We have the following control on the kernel (see Lemma 5.5 of \cite{dela:meno:10} for a proof).
\begin{LEMME}
\label{lem:F:050509:2}
Let $T_0>0$ be fixed. There exists 
a constant $C_{\ref{lem:F:050509:2}} >0$,
depending on ({\bf A}-$\eta $) and $T_0$
 such that, for all
$t \in [0,T),\ T\le T_0$ and
$\z,\y \in \R^{nd}$, 
\begin{equation*}
|H(t,T,\z,\by)|
\leq C_{\ref{lem:F:050509:2}}
(T-t)^{\frac \eta 2 -1}
g_{C_{\ref{lem:F:050509:2}},T-t}
\bigl(\z - \btheta_{t,T}(\y) \bigr).
\end{equation*}
\end{LEMME}
We conclude this section with a technical Lemma whose proof is postponed to Appendix \ref{APP}.
\begin{LEMME}
\label{LEMME_CV_DIR}
Let  $h $ be a $C^0(\R^{nd},\R)$ function. Define for all $(s,\x)\in [0,T]\times \R^{nd}$,
\begin{equation*}
\forall \varepsilon >0,\ \bTheta^\varepsilon(s,\x):=\bint{\R^{nd}}^{} d\y h(s,\y) \tilde p^{s+\varepsilon,\y}(s,s+\varepsilon,\x,\y).
\end{equation*} 
Then
$\bTheta^\varepsilon(s,x) $ converges boundedly and pointwise to $h(s,\x) $ when $\varepsilon \rightarrow 0 $.
\end{LEMME}

\mysection{Proof of Theorem \ref{THM_M}}
Suppose we are given two solutions $\P_1,\P_2$ of the martingale problem associated to $(L_t)_{t\in [s,T]}$ starting in $\x$ at time $s$. W.l.o.g. we can suppose here that $T\le 1 $.
Define for a bounded Borel function $f:[0,T]\times \R^{nd}\rightarrow \R$:
$$S^i f:=\E_i\left[\int_{s}^{T} dt f(t,\X_t) \right],\ i\in \leftB 1, 2\rightB, $$
where $(\X_t)_{t\in [s,T]}$ stands here for the canonical process associated to $(\P_i)_{i\in \leftB 1,2\rightB} $. Let us specify (as indicated in \cite{bass:perk:09}) that $Sf$ is only a linear functional and not a function since $\P_i$ does not need to come from a Markov process. 
Let us now introduce 
$$S^\Delta f:=S^1 f-S^2f,\ \Theta:=\sup_{\|f\|_{\infty}\le 1} |S^\Delta f|.$$
Clearly $\Theta\le T-s$.

If $f\in C_0^2([0,T)\times \R^{nd},\R)$ then by definition of the martingale problem we have:
\begin{equation}
\label{DEF_PM}
f(s,\x)+\E_i\left[ \int_s^T dt (\partial_t +L_t) f(t,\X_t)\right]=0,\ i\in \leftB 1,2\rightB.
\end{equation}

For a fixed point $\y\in \R^{nd}$ and $\varepsilon \ge 0 $, introduce $\forall f\in C_0^{1,2}([0,T)\times \R^{nd},\R)  $ the Green function
\begin{equation}
\label{GREEN}
\forall (s,\x) \in [0,T)\times \R^{nd},\ G^{\varepsilon,\y} f(s,\x):=\int_s^T dt \int_{\R^{nd}}^{}d\z  \tilde p^{t+\varepsilon,\y}(s,t,\x,\z) f(t,\z).
\end{equation}
We insist that in the above equation $\tilde p^{t+\varepsilon,\y}(s,t,\x,\z) $ stands for the density at time $t$ and point $\z$ of the process $\tilde \X^{t+\varepsilon,\y}$ defined in \eqref{eq:F:linear:eq:upper} starting from $\x$ at time $s$ with coefficients depending on the backward transport of the freezing point $\y$ by $(\btheta_{u,t+\varepsilon})_{u\in [s,t]} $. In particular, the parameter $\varepsilon $ can be equal to $0 $ in the previous definition.

One easily checks that
\begin{equation}
\label{EQ_FK_GEL}
(\partial_s+\tilde L_{s,\x}^{t+\varepsilon,\y})\tilde p^{t+\varepsilon,\y}(s,t,\x,\z)=0, \forall (s,\x,\z) \in [0,t)\times (\R^{nd})^2,\ \tilde p^{t+\varepsilon,\y}(s,t,\x,.)\underset{t\downarrow s}{\longrightarrow} \delta_{\x}(.). 
\end{equation}
Introducing for all $f\in C_0 ^{1,2}([0,T)\times \R^{nd},\R),\ (s,x)\in [0,T)\times \R^{nd} $,
$$M_{s,\x}^{\varepsilon,\y} f(s,\x):=\int_{s}^{T}dt \int_{\R^{nd}}^{} d\z \tilde L_{s,\x}^{t+\varepsilon,\y} \tilde p^{t+\varepsilon,\y}(s,t,\x,\z) f(t,\z),$$
we get from equations \eqref{GREEN}, \eqref{EQ_FK_GEL}  
\begin{equation}
\label{EDP_GREEN}
(\partial_s G^{\varepsilon,\y} f+M_{s,\x}^{\varepsilon,\y}f) (s,\x)=-f(s,\x), \forall (s,\x)\in [0,T)\times \R^{nd}.
\end{equation}
Define now for a smooth function $h\in C_0^{1,2}([0,T)\times \R^{nd},\R) $ and for all $ (s,\x) \in [0,T)\times \R^{nd}$:
\begin{eqnarray*}
\Phi^{\varepsilon,\y}(s,\x)&:=&\tilde p^{s+\varepsilon,\y}(s,s+\varepsilon,\x,\y) h(s,\y),\\
\Psi_\varepsilon(s,\x)&:=&\bint{\R^{nd}}^{} d\y G^{\varepsilon,\y} (\Phi^{\varepsilon,\y})(s,\x).
\end{eqnarray*}
Observe that \eqref{GREEN} yields:
\begin{eqnarray}
\Psi_\varepsilon(s,\x)&:=&\bint{\R^{nd}}^{}d \y \bint{s}^{T} dt \bint{\R^{nd}}^{} d\z \tilde p^{t+\varepsilon,\y}(s,t,\x,\z)\Phi^{\varepsilon,\y}(t,\z)\nonumber \\
&=&\bint{\R^{nd}}^{}d \y \bint{s}^{T} dt \bint{\R^{nd}}^{} d\z \tilde p^{t+\varepsilon,\y}(s,t,\x,\z)\tilde p^{t+\varepsilon,\y}(t,  t+\varepsilon,\z,\y)h(t,\y)\nonumber \\
&=&\bint{\R^{nd}}^{}d\y  \bint{s}^{T} dt  \tilde p^{t+\varepsilon,\y}(s,t+\varepsilon,\x,\y) h(t,\y), \label{EXPR_PHI_EPS}
\end{eqnarray}
exploiting the semigroup property of the frozen density $\tilde p^{t+\varepsilon,\y} $ for the last inequality.
Write now for all $(s,x)\in [0,T)\times \R^{nd}$,
\begin{eqnarray*}
(\partial_s+L_s) \Psi_\varepsilon (s,x)&=&\bint{\R^{nd}}^{}  d\y  (\partial_s+L_s)(G^{\varepsilon,\y}\Phi^{\varepsilon,\y}   )(s,\x)\\
&=&\bint{\R^{nd}}^{}  d\y  (\partial_s G^{\varepsilon,\y}\Phi^{\varepsilon,\y}+M_{s,\x}^{\varepsilon,\y} \Phi^{\varepsilon,\y})(s,\x)+\bint{\R^{nd}}^{}d\y (L_s G^{\varepsilon,\y}\Phi^{\varepsilon,\y}-M_{s,\x}^{\varepsilon,\y}\Phi^{\varepsilon,\y})(s,\x)\\
&=&-\int_{\R^{nd}}^{}d\y \Phi^{\varepsilon,\y}(s,\x)+\bint{\R^{nd}}^{}d\y (L_s G^{\varepsilon,\y}\Phi^{\varepsilon,\y}-M_{s,\x}^{\varepsilon,\y}\Phi^{\varepsilon,\y})(s,\x):=I_1^\varepsilon+I_2^\varepsilon,
\end{eqnarray*}
exploiting \eqref{EDP_GREEN} for the last but one identity. Now Lemma \ref{LEMME_CV_DIR} gives 
%\begin{eqnarray*}
%\label{CT_I1}
$I_1^\varepsilon \underset{\varepsilon\rightarrow 0}{\longrightarrow} -h(s,\x)$. 
%\end{eqnarray*}
On the other hand using the notations of Section \ref{SEC_FRO}, we derive from \eqref{EXPR_PHI_EPS} that the term $I_2^\varepsilon $ writes:
\begin{eqnarray*}
I_2^\varepsilon &=&\int_s^T dt \bint{\R^{nd}}^{} d\y (L_s-\tilde L_{s,\x}^{t+\varepsilon,\y}) \tilde p^{t+\varepsilon,\y}(s,t+\varepsilon,\x,\y) h(t,\y)\\
                         &
                         =
                         & 
                         \int_s^T dt \bint{\R^{nd}}^{} d\y H(s,t+\varepsilon,\x,\y) h(t,\y).
\end{eqnarray*}
Thus, from Lemma \ref{lem:F:050509:2},
\begin{eqnarray}
|I_2^\varepsilon|&\le& C_{\ref{lem:F:050509:2}}  \|h\|_{\infty} \bint{s}^T dt (t+\varepsilon-s)^{-1+\eta/2} \bint{\R^{nd}}^{}d\y g_{C_{\ref{lem:F:050509:2}},t+\varepsilon-s}(\x-\btheta_{s,t+\varepsilon} (\y))\nonumber \\
&\le & \|h\|_\infty C_{\ref{CT_I2}}\bint{s}^T  dt (t+\varepsilon-s)^{-1+\eta/2} \bint{\R^{nd}}^{}d\y g_{C_{\ref{CT_I2}},t+\varepsilon-s}(\btheta_{t+\varepsilon,s}(\x) -\y) \le C_{\ref{CT_I2}} \{(T-s)\vee \varepsilon \}^{\eta/2} \|h\|_\infty,\nonumber\\
\label{CT_I2}
\end{eqnarray}
using the bi-Lipschitz property of the flow for the last but one inequality and up to a modification of $C_{\ref{CT_I2}} $ in the last one. Anyhow, the constant $C_{\ref{CT_I2}} $ only depends on known parameters in \A{A-$\eta$}. Thus for $T$ and $\varepsilon  $ sufficiently small we have from \eqref{CT_I2}
\begin{eqnarray}
\label{CT_I2_DEF}
|I_2^\varepsilon|\le \frac 12 \|h\|_\infty.
\end{eqnarray}
Now, equation \eqref{DEF_PM} and the above definition of $S^\Delta $ yield:
\begin{eqnarray*}
S^\Delta((\partial_s+L_s) \Psi_\varepsilon)=0 \Rightarrow |S^\Delta I_1^\varepsilon|=|S^\Delta I_2^\varepsilon|.
\end{eqnarray*}
From the bounded convergence part of Lemma \ref{LEMME_CV_DIR} and \eqref{CT_I2_DEF}, we have:
\begin{eqnarray*}
|S^\Delta h|=\lim_{\varepsilon\rightarrow 0} |S^\Delta I_1^\varepsilon|=\lim_{\varepsilon\rightarrow 0} |S^\Delta I_2^\varepsilon|\le \Theta \limsup_{\varepsilon \rightarrow 0} |I_2^\varepsilon|\le \Theta \frac{\|h\|_\infty}2.
\end{eqnarray*}
By a monotone class argument, the previous inequality remains valid for bounded measurable functions $h$ compactly supported in $[0,T)\times \R^{nd} $. Taking the supremum over $\|h\|_\infty\le 1 $, we obtain $\Theta \le \frac 12 \Theta$ which gives $\Theta=0 $ since $\Theta<+\infty $. Hence, %for all  $g\in C_0^0([0,T)\times \R^{nd}, \R^{nd}),\ 
$\E_1\left[\int_{s}^T dt h(t,\X_t) \right]=\E_2\left[\int_{s}^T dt h(t,\X_t) \right] $ which proves the result on the interval $[0,T] $. Regular conditional probabilities then allow to extend the result on $\R^+ $, see Chapter 6.2 of \cite{stro:vara:79} for details.

\appendix

\mysection{Proof  of Lemma \ref{LEMME_CV_DIR}}
\label{APP}
Let us denote by $\tilde \K_{s,s+\varepsilon}^{s+\varepsilon,\y} $ the covariance matrix associated to equation \eqref{eq:F:linear:eq:upper} for the process $\tilde \X^{s+\varepsilon,\y} $ (starting from $\x$ at $s$) at time $s+\varepsilon $ and by $\K_{s,s+\varepsilon}^{s,\x}$ the covariance matrix associated to the linear diffusion with dynamics:
\begin{equation}
\label{BAR_X}
d\bar \X_t=D\gF(t,\btheta_{t,s}( \x))\bar \X_t dt +B\sigma(t,\btheta_{t,s}(\x)) dW_t,  t\ge s,
\end{equation}
that is $\K_{s,s+\varepsilon}^{s,\x}=\int_s^{s+\varepsilon} du \gR^{s,\x}(s+\varepsilon,u) Ba^*(u,\btheta_{u,s}(\x)) B^* \gR^{s,\x}(s+\varepsilon,u)^* du $ where $ \gR^{s,\x}$ stands for the resolvent 
associated to the linear part of \eqref{BAR_X}.

Under \A{A-$\eta$}, the matrices $\tilde \K_{s,s+\varepsilon}^{s+\varepsilon,\y},\ \K_{s,s+\varepsilon}^{s,\x} $ admit a good scaling property in the sense of the previous Definition \ref{GSP}, i.e. 
\begin{eqnarray}
\label{ND_COND}
\exists C:=C(\A{A{\rm {-}}\eta})\ge 1,\ \forall \bxi \in \R^{nd},&&C^{-1}\varepsilon^{-1}|\T_{\varepsilon} \bxi|^2\le  \langle \tilde \K_{s,s+\varepsilon}^{s+\varepsilon,\y} \bxi, \bxi\rangle  \le C\varepsilon^{-1}|\T_{\varepsilon} \bxi|^2, \nonumber\\ 
&&   C^{-1}\varepsilon^{-1}|\T_{\varepsilon} \bxi|^2\le  \langle \K_{s,s+\varepsilon}^{s,\x} \bxi, \bxi\rangle  \le C\varepsilon^{-1}|\T_{\varepsilon} \bxi|^2. 
\end{eqnarray}

We introduce the following decomposition:
\begin{eqnarray}
\bTheta^\varepsilon(s,\x):=\bint{\R^{nd}}^{} \frac{d\y}{(2\pi)^{nd/2}} h(s,\y)\exp\left (-\frac 12 \langle (\tilde \K_{s,s+\varepsilon}^{s+\varepsilon,\y})^{-1}(\tilde \btheta_{s+\varepsilon,s}^{s+\varepsilon,\y} (\x)-\y), \tilde \btheta_{s+\varepsilon,s}^{s+\varepsilon,\y} (\x)-\y\rangle  \right)\nonumber\\
\times \left\{\frac{1}{\det(\tilde \K_{s,s+\varepsilon}^{s+\varepsilon,\y})^{1/2}}-\frac{1}{ \det(\K_{s,s+\varepsilon}^{s ,\x})^{1/2}} \right\}\nonumber \\
+\bint{\R^{nd}}^{}\frac{d\y}{(2\pi)^{nd/2}} \frac{h(s,\y)}{\det(\K_{s,s+\varepsilon}^{s,\x})^{1/2}}
\left[ \exp\left (-\frac 12 \langle (\tilde \K_{s,s+\varepsilon}^{s+\varepsilon,\y})^{-1}(\tilde \btheta_{s+\varepsilon,s}^{s+\varepsilon,\y} (\x)-\y), \tilde \btheta_{s+\varepsilon,s}^{s+\varepsilon,\y} (\x)-\y\rangle  \right)-\right. \nonumber \\
\left. \exp\left (-\frac 12 \langle  (\K_{s,s+\varepsilon}^{s,\x})^{-1}(\tilde \btheta_{s+\varepsilon,s}^{s+\varepsilon,\y} (\x)-\y), \tilde \btheta_{s+\varepsilon,s}^{s+\varepsilon,\y} (\x)-\y\rangle  \right)\right]\nonumber\\
+\bint{\R^{nd}}^{}\frac{d\y}{(2\pi)^{nd/2}} \frac{h(s,\y)}{\det(\K_{s,s+\varepsilon}^{s,\x})^{1/2}}
\left[ \exp\left (-\frac 12 \langle  (\K_{s,s+\varepsilon}^{s,\x})^{-1}(\tilde \btheta_{s+\varepsilon,s}^{s+\varepsilon,\y} (\x)-\y), \tilde \btheta_{s+\varepsilon,s}^{s+\varepsilon,\y} (\x)-\y\rangle  \right)-\right. \nonumber\\
\left. \exp\left (-\frac 12 \langle  (\K_{s,s+\varepsilon}^{s,\x})^{-1}(\btheta_{s+\varepsilon,s} (\x)-\y), \btheta_{s+\varepsilon,s}(\x)-\y\rangle  \right)\right]\nonumber\\
+\bint{\R^{nd}}^{}\frac{d\y}{(2\pi)^{nd/2}} \frac{h(s,\y)}{\det(\K_{s,s+\varepsilon}^{s,\x})^{1/2}}\exp\left (-\frac 12 \langle  (\K_{s,s+\varepsilon}^{s,\x})^{-1}(\btheta_{s+\varepsilon,s} (\x)-\y), \btheta_{s+\varepsilon,s}(\x)-\y\rangle  \right):=\bsum{i=1}^{4}\bTheta_i^\varepsilon(s,\x).\label{decomp_bTHETA}
\end{eqnarray}

Let us now control the $(\bTheta_i^\varepsilon(s,\x))_{i\in \leftB 1,4\rightB} $. Set $\tilde \y:=(\K_{s,s+\varepsilon }^{s,\x} )^{-1/2}(\btheta_{s+\varepsilon,s}(\x)-\y)$ 
where $(\K_{s,s+\varepsilon}^{s,\x})^{-1/2}$ denotes the upper triangular matrix obtained through the Cholesky factorization, i.e. $ (\K_{s,s+\varepsilon}^{s,\x})^{-1}=((\K_{s,s+\varepsilon}^{s,\x})^{-1/2})^*(\K_{s,s+\varepsilon}^{s,\x})^{-1/2} $. We get from the bounded convergence theorem:
\begin{eqnarray}
\label{CTR_TH_4}
\bTheta_4^\varepsilon(s,\x)=\bint{\R^{nd}}^{} \frac{d\tilde \y}{(2\pi)^{nd/2}} h(s,-(\K_{s,s+\varepsilon}^{s,\x})^{1/2}\tilde \y+\btheta_{s+\varepsilon,s}(\x)) \exp\left(-\frac{|\tilde \y|^2}{2} \right) \underset{\varepsilon \rightarrow 0}{\longrightarrow} h(s,\x). 
\end{eqnarray}

For $\bTheta_3^\varepsilon(s,\x) $, we first observe, from Lemma \ref{eq:F:240409:2}, the good scaling property \eqref{ND_COND} %(deriving from assumption \A{A}) 
and the bi-Lipschitz property of the flow $\btheta$, that there exists $C_1:=C_1($\A{A-$\eta $}$)\ge 1$ s.t. 
\begin{eqnarray}
\label{COND_EQUIV}
C_1^{-1}\varepsilon^{1/2}|\T_\varepsilon^{-1}(\tilde \btheta_{s+\varepsilon,s}^{s+\varepsilon,\y}(\x)-\y)|\le | (\K_{s,s+\varepsilon}^{s,\x})^{-1/2}(  \tilde \btheta_{s+\varepsilon,s}^{s+\varepsilon,\y}(\x)-\y)| \le C_1\varepsilon^{1/2}|\T_\varepsilon^{-1}(\tilde \btheta_{s+\varepsilon,s}^{s+\varepsilon,\y}(\x)-\y)|,\nonumber\\
C_1^{-1}\varepsilon^{1/2}|\T_\varepsilon^{-1}( \btheta_{s+\varepsilon,s}(\x)-\y)|\le | (\K_{s,s+\varepsilon}^{s,\x})^{-1/2}  (\btheta_{s+\varepsilon,s} (\x)-\y)| \le C_1\varepsilon^{1/2}|\T_\varepsilon^{-1}(\btheta_{s+\varepsilon,s}(\x)-\y)|,\nonumber\\
C_1^{-1} \varepsilon^{1/2} | \T_\varepsilon^{-1}( \btheta_{s+\varepsilon,s}(\x)-\y) |\le  \varepsilon^{1/2} | \T_\varepsilon^{-1}(\tilde \btheta_{s+\varepsilon,s}^{s+\varepsilon,\y}(\x)-\y) |
 \le C_1\varepsilon^{1/2} | \T_\varepsilon^{-1}( \btheta_{s+\varepsilon,s}(\x)-\y) |.
\end{eqnarray}

%%%%%% Je ne pense pas que cela serve en fait 
%For a given $N\ \in \R $, we thus get $B_N:=\{\y\in \R^{nd}: \varepsilon^{1/2}|\T_\varepsilon^{-1}(\btheta_{s+\varepsilon,s}(\x)-\y)|\ge N \}\subset \{\y \in \R^{nd}: \varepsilon^{1/2}|\T_\varepsilon^{-1}(\tilde \btheta_{s+\varepsilon,s}^{s+\varepsilon,\y}(\x)-\y)|\ge N C_1^{-1}\}$. Thus,  writing $\bTheta_3^\varepsilon(s,\x):=\int_{\R^{nd}}^{} d\y (G_1^\varepsilon(s,\x,\y)-G_2^\varepsilon(s,\x,\y)) $ we then derive that there exists $C_2:=C_2(\A{A})$ s.t.
%\begin{eqnarray}
%\label{CTR_TH_4_BOUTS}
%\bTheta_{32}^\varepsilon(s,\x):=\bint{B_N}^{} d\y (|G_1^\varepsilon|+|G_2^\varepsilon|)(s,\x,\y) \le C_2\exp(-C_2^{-1}N^2).
%\end{eqnarray}
Write now,
\begin{eqnarray}
%\bTheta_3^\varepsilon(s,\x)&:=&\int_{B_N^C}^{} d\y(G_1^\varepsilon-G_2^\varepsilon)(s,\x,\y) +\bTheta_{32}(s,\x):=(\bTheta_{31}+\bTheta_{32})(s,\x).\\
|\bTheta_{3}^\varepsilon(s,\x)|&\le & \|h\|_\infty \bint{\R^{nd}}^{} \frac{d\y}{(2\pi)^{nd/2}\det(\K_{s,s+\varepsilon}^{s,\x})^{1/2}} \int_0^1 d\delta |(\varphi_{s,\x,\y}^\varepsilon)'(\delta)|,\ \forall  \delta \in [0,1],\label{THE_FORM_R3} \\
\varphi_{s,\x,\y}^\varepsilon(\delta)&=&\exp\left(-\frac{1}{2}\left\{\langle (\K_{s,s+\varepsilon}^{s,\x})^{-1} (\tilde \btheta_{s+\varepsilon,s}^{s+\varepsilon,\y}(\x)-\y),\tilde \btheta_{s+\varepsilon,s}^{s+\varepsilon,\y}(\x)-\y  \rangle +\right. \right. \nonumber\\
&& \left. \left. \delta \left[\langle (\K_{s,s+\varepsilon}^{s,\x})^{-1} ( \btheta_{s+\varepsilon,s}(\x)-\y), \btheta_{s+\varepsilon,s}(\x)-\y  \rangle-\langle (\K_{s,s+\varepsilon}^{s,\x})^{-1} (\tilde \btheta_{s+\varepsilon,s}^{s+\varepsilon,\y}(\x)-\y),\tilde \btheta_{s+\varepsilon,s}^{s+\varepsilon,\y}(\x)-\y  \rangle \right]\right\}\right),\nonumber\\
|(\varphi_{s,\x,\y}^\varepsilon)'(\delta)|&\le& \left | (\K_{s,s+\varepsilon}^{s,\x})^{-1/2}\left\{(\btheta_{s+\varepsilon,s}(\x) -\y)+(\tilde \btheta_{s+\varepsilon,s}^{s+\varepsilon,\y}(\x) -\y) \right\}\right|\nonumber\\
&& \times \left|(\K_{s,s+\varepsilon}^{s,\x})^{-1/2}\left\{(\btheta_{s+\varepsilon,s}(\x) -\y)-(\tilde \btheta_{s+\varepsilon,s}^{s+\varepsilon,\y}(\x) -\y) \right\} \right|  \varphi_{s,\x,\y}^\varepsilon(\delta),\nonumber
\end{eqnarray} 
using the Cauchy-Schwarz inequality for the last assertion. Equations \eqref{COND_EQUIV} now yield that there exists $C_2:=C_2($\A{A-$\eta$}$)\ge 1$ s.t.:
\begin{eqnarray}
\label{CTR_DER_PHI}
|(\varphi_{s,\x,\y}^\varepsilon)'(\delta)| &\le& C_2 \varepsilon^{1/2} |\T_\varepsilon^{-1}\left(\tilde \btheta_{s+\varepsilon,s}^{s+\varepsilon,\y}(\x)-\btheta_{s+\varepsilon,s}(\x)\right)| \exp\left(-C_2^{-1}\varepsilon |\T_\varepsilon^{-1}(\btheta_{s+\varepsilon,s}(\x)-\y)|^2 \right)\nonumber\\
&:=& C_2 |D_\varepsilon| \exp\left(-C_2^{-1}\varepsilon |\T_\varepsilon^{-1}(\btheta_{s+\varepsilon,s}(\x)-\y)|^2 \right).
\end{eqnarray}
Let us now recall the differential dynamics of $ \btheta_{s+\varepsilon,s}(\x),\ \tilde \btheta_{s+\varepsilon,s}^{s+\varepsilon,\y}(\x)$, that is $\btheta_{s+\varepsilon,s}(\x)=\x+\int_s^{s+\varepsilon} du \gF(u,\btheta_{u,s}(\x)) $ and $\tilde \btheta _{s+\varepsilon,s}^{s+\varepsilon,\y}(\x)=\x+\int_{s}^{s+\varepsilon} du\{\gF(u,\btheta_{u,s+\varepsilon}(\y))+ D\gF(u,\btheta_{u,s+\varepsilon}(\y))(\tilde \btheta_{u,s}^{s+\varepsilon,\y}(\x)-\btheta_{u,s+\varepsilon}(\y))\}$. Set %for all $\z\in \R^{nd}$:
$ \gF^{s+\varepsilon,\y}(u,\z):=$
$(\gF_1(u,\btheta_{u,s+\varepsilon}(\y) ),\gF_2(u,\z_1, (\btheta_{u,s+\varepsilon}(\y))^{2,n}),\gF_3(u,\z_2,(\btheta_{u,s+\varepsilon}(\y))^{3,n})),\cdots,$ 
$\gF_n(u,\z_{n-1},(\btheta_{u,s+\varepsilon}(\y))_n ) )^*,\ \forall  (u,\z)\in [s,s+\varepsilon]\times \R^{nd}$. Observe in particular that $\gF^{s+\varepsilon,\y}(u,\btheta_{u,s+\varepsilon}(\y))=\gF(u,\btheta_{u,s+\varepsilon}(\y)) $. 
We get:
\begin{eqnarray}
D_\varepsilon:=\varepsilon^{1/2} \T_{\varepsilon}^{-1}\left\{\btheta_{s+\varepsilon,s}(\x)-\tilde \btheta_{s+\varepsilon,s}^{s+\varepsilon,\y}(\x) \right\}=\nonumber\\
\varepsilon^{1/2} \T_{\varepsilon}^{-1}\left\{ \int_{s}^{s+\varepsilon} du \biggl[\biggl( \gF(u,\btheta_{u,s}(\x))-\gF^{s+\varepsilon,\y}(u,\btheta_{u,s}(\x)) \biggr)+\biggl(D\gF(u,\btheta_{u,s+\varepsilon}(\y))(\btheta_{u,s}(\x)-\tilde \btheta_{u,s}^{s+\varepsilon,\y}(\x) )\biggr) \right. \nonumber \\
\left. + \biggl(\bint{0}^{1} d\delta \left(  D\gF^{s+\varepsilon,\y}(u,\btheta_{u,s+\varepsilon}(\y)+\delta (\btheta_{u,s}(\x)-\btheta_{u,s+\varepsilon}(\y) ) ) -D\gF^{s+\varepsilon,\y}(u,\btheta_{u,s+\varepsilon}(\y)) \right) (\btheta_{u,s}(\x)-\btheta_{u,s+\varepsilon}(\y) )\biggr) \biggr] \right\}\nonumber \\
:=D_\varepsilon^1+D_\varepsilon^2+D_\varepsilon^3,\nonumber\\  
\label{D_EPS}
\end{eqnarray}
where for $(u,\z)\in [s,s+\varepsilon]\times \R^{nd}, D\gF^{s+\varepsilon,\y}(u,\z) $ is the $(nd)\times (nd)$ matrix with only non zero $d\times d $ matrix entries $(D\gF^{s+\varepsilon}(u,\z))_{j,j-1}:=D_{\x_{j-1}}\gF_j(u,\z_{j-1},\btheta_{u,s+\varepsilon}(\y)^{j,n}) ,\  j\in \leftB 2,n\rightB $, so that in particular $D\gF^{s+\varepsilon,\y}(u,\btheta_{u,s+\varepsilon}(\y))=D\gF(u,\btheta_{u,s+\varepsilon}(\y)) $.

The structure of the ``partial gradient" $D\gF^{s+\varepsilon,\y} $ associated to the $\eta $-H\"older continuity of the mapping $ \x_{j-1} \in \R^d \mapsto D_{\x_{j-1}}\gF_j(\x_{j-1},\x^{j,n} ),\forall \x^{j,n}\in \R^{(n-j+1)d}  $ yield that there exists $C_3:=C_3($\A{A-$\eta$} $)$ s.t. for all $j\in \leftB 2,d \rightB $:
\begin{eqnarray}
|(D_\varepsilon^3)_j| &\le & C_3 \varepsilon^{1/2-j}\int_{s}^{s+\varepsilon} du|(\btheta_{u,s}(\x)- \btheta_{u,s+\varepsilon}(\y))_{j-1}|^{1+\eta}\nonumber \\
&\le & C_3 \varepsilon^{-1}\int_s^{s+\varepsilon}du (\sum_{k=2}^n  \varepsilon^{1/2-(k-1)}|(\btheta_{u,s}(\x)- \btheta_{u,s+\varepsilon}(\y))_{k-1}| )^{1+\eta}\varepsilon^{((j-1)-1/2)\eta}\nonumber \\
&\le & C_3 \varepsilon^{-1+\eta(j-3/2)}\int_{s}^{s+\varepsilon} du (\varepsilon^{1/2}|\T_\varepsilon^{-1}(\btheta_{u,s}(\x)- \btheta_{u,s+\varepsilon} (\y))|)^{1+\eta}\nonumber \\
&\le & C_3 \varepsilon^{-1+\eta(j-3/2)}\int_{s}^{s+\varepsilon} du (\varepsilon^{1/2}|\T_\varepsilon^{-1}(\btheta_{s+\varepsilon,s}(\x)-\y)|)^{1+\eta}\nonumber \\
&\le & C_3 \varepsilon^{\eta(j-3/2)}(\varepsilon^{1/2}|\T_\varepsilon^{-1}(\btheta_{s+\varepsilon,s}(\x)-\y)|)^{1+\eta}, \label{CTR_D3}
\end{eqnarray}
up to a modification of $C_3$ and using  the bi-Lipschitz property of the flow $\btheta$ for the last but one inequality (see the end of the proof of Proposition 5.1 in \cite{dela:meno:10} for details).

On the other hand, the term $D_\varepsilon^1$ can be seen as a remainder w.r.t. the characteristic time scales. Precisely, there exists $C_4:=C_4($\A{A-$\eta$}$)$ s.t. for all $j\in \leftB 1,n \rightB $:
\begin{eqnarray}
|(D_\varepsilon^1)_j|&\le&C_4\varepsilon^{1/2-j}\int_s^{s+\varepsilon}du \sum_{k=j}^n|(\btheta_{u,s}(\x)-\btheta_{u,s+\varepsilon}(\y))_k|\nonumber \\
                                 &\le & C_4\int_{s}^{s+\varepsilon} du \varepsilon^{1/2}|\T_\varepsilon^{-1}(\btheta_{u,s}(\x)-\btheta_{u,s+\varepsilon}(\y))| \nonumber \\
                                 &\le & C_4 \varepsilon (\varepsilon^{1/2}|\T_\varepsilon^{-1}(\btheta_{s+\varepsilon,s}(\x)-\y)|)\label{CTR_D1}
\end{eqnarray}
using once again the  bi-Lipschitz property of the flow $\btheta $ for the last inequality. Recall now that $D_2^\varepsilon$ is the linear part of equation \eqref{D_EPS}, i.e. it can be rewritten 
\begin{eqnarray*}
D_2^\varepsilon&=&\int_s^{s+\varepsilon} du \left\{\varepsilon^{1/2}\T_\varepsilon^{-1} D\gF(u,\btheta_{u,s+\varepsilon}(\y)) (u-s)^{-1/2}\T_{u-s}\right\}\left((u-s)^{1/2}\T_{u-s}^{-1}(\btheta_{u,s}(\x)-\tilde \btheta_{u,s}^{s+\varepsilon,\y}(\x)) \right)\\ 
&:=&\int_s^{s+\varepsilon} du \alpha_\varepsilon^\y(u,s)\left((u-s)^{1/2}\T_{u-s}^{-1} (\btheta_{u,s}(\x)-\tilde \btheta_{u,s}^{s+\varepsilon,\y}(\x))\right),
\end{eqnarray*}
where there exists a constant $\tilde C:=\tilde C($\A{A-$\eta $}) independent of $\varepsilon $ s.t. $\int_{s}^{s+\varepsilon} du|\alpha_\varepsilon^\y(u,s)|\le \tilde C $. %can be bounded independently of $\varepsilon $. 
From \eqref{CTR_D1}, \eqref{CTR_D3}, \eqref{D_EPS} and Gronwall's Lemma  we derive 
$$\exists C_5:=C_5(\A{A-\eta}),\ |D_\varepsilon|\le C_5 \varepsilon^{\eta/2}((\varepsilon^{1/2}|\T_\varepsilon^{-1}(\btheta_{s+\varepsilon,s}(\x)-\y)|)^{1+\eta}+1).$$
Plugging this estimate into \eqref{CTR_DER_PHI}, we then get  from  \eqref{THE_FORM_R3}, using as well the good scaling property \eqref{ND_COND}, that there exists $C_6:=C_6($\A{A-$\eta$}$)$,
\begin{eqnarray}
|\bTheta_{3}^\varepsilon(s,\x)|&\le& C_6 \varepsilon^{\eta/2}\int_{\R^{nd}}^{} \frac{d\y}{\varepsilon^{n^2d/2}} ((\varepsilon^{1/2}|\T_\varepsilon^{-1}(\btheta_{s+\varepsilon,s}(\x)- \y)|)^{1+\eta}+1)\exp(-C_2^{-1}\varepsilon|\T_\varepsilon^{-1}(\btheta_{s+\varepsilon,s}(\x)- \y)|^2)\nonumber\\
&\le& C_6 \varepsilon^{\eta/2},
\label{CTR_R3}
\end{eqnarray}
up to a modification of $C_6$ in the last inequality. %%%%% Changer tout cela 

Let us consider now $\bTheta_2^\varepsilon(s,\x) $. Write:
\begin{eqnarray}
|\bTheta_2^\varepsilon(s,\x)|&\le &C_7 \bint{\R^{nd}}^{}\frac{d\y}{(2\pi)^{nd/2}} \frac{h(s,\y)}{\det(\K_{s,s+\varepsilon}^{s,\x})^{1/2}}\int_0^1 d\delta |(\psi_{s,\x,\y}^\varepsilon)'(\delta)|,\ C_7:=C_7(\A{A-\eta}),\ \forall \delta \in [0,1], \label{DEF_BTHETA_2}\\
\psi_{s,\x,\y}^\varepsilon(\delta)&=&\exp\left(-\frac 12 \biggl\{ \langle (\K_{s,s+\varepsilon}^{s,\x})^{-1} (\tilde \btheta_{s+\varepsilon,s}^{s+\varepsilon,\y}(\x)-\y), \tilde \btheta_{s+\varepsilon,s}^{s+\varepsilon,\y}(\x)-\y\rangle    \right.\nonumber\\
&+&\left.\delta \biggl[ \langle (\tilde \K_{s,s+\varepsilon}^{s+\varepsilon,\y})^{-1} (\tilde \btheta_{s+\varepsilon,s}^{s+\varepsilon,\y}(\x)-\y), \tilde \btheta_{s+\varepsilon,s}^{s+\varepsilon,\y}(\x)-\y\rangle  -\langle (\K_{s,s+\varepsilon}^{s,\x})^{-1} (\tilde \btheta_{s+\varepsilon,s}^{s+\varepsilon,\y}(\x)-\y), \tilde \btheta_{s+\varepsilon,s}^{s+\varepsilon,\y}(\x)-\y\rangle\biggl] \biggr\}   \right),\nonumber \\
|(\psi_{s,\x,\y}^\varepsilon)'(\delta)|&\le&  \left |  \langle   ( (\tilde \K_{s,s+\varepsilon}^{s+\varepsilon,\y})^{-1}-(\K_{s,s+\varepsilon}^{s,\x})^{-1})(\tilde \btheta_{s+\varepsilon,s}^{s+\varepsilon,\y}(\x)-\y), \tilde \btheta_{s+\varepsilon,s}^{s+\varepsilon,\y}(\x)-\y \rangle \right | %\\
%&&\times \biggl(|(\tilde \K_{s+\varepsilon,s}^{s+\varepsilon,\y})^{-1/2}(\tilde \btheta_{s+\varepsilon,s}^{s+\varepsilon,\y}(\x)-\y)|+|(\K_{s+\varepsilon,s}^{s,\x})^{-1/2}(\tilde \btheta_{s+\varepsilon,s}^{s+\varepsilon,\y}(\x)-\y)|\biggr)
\psi_{s,\x,\y}^\varepsilon(\delta)\nonumber.
\end{eqnarray}
Equations \eqref{COND_EQUIV} (that according to \eqref{ND_COND} hold for $\tilde \K_{s+\varepsilon,s}^{s,\y} $ as well) yield:
\begin{eqnarray}
|(\psi_{s,\x,\y}^\varepsilon)'(\delta)|&\le& C \left |  \langle   ( (\tilde \K_{s,s+\varepsilon}^{s+\varepsilon,\y})^{-1}-(\K_{s,s+\varepsilon}^{s,\x})^{-1})(\tilde \btheta_{s+\varepsilon,s}^{s+\varepsilon,\y}(\x)-\y), \tilde \btheta_{s+\varepsilon,s}^{s+\varepsilon,\y}(\x)-\y \rangle \right |\exp(-C\varepsilon |\T_{\varepsilon}^{-1}(\btheta_{s+\varepsilon,s}(\x)-\y)|^2)\nonumber \\
&:=&C |Q_\varepsilon| \exp(-C\varepsilon |\T_{\varepsilon}^{-1}( \btheta_{s+\varepsilon,s}(\x)-\y)|^2),\label{CTR_PSI_DER}
\end{eqnarray}
for $C:=C($\A{A-$\eta$}$)$. From the scaling Lemma 3.6  in \cite{dela:meno:10}, we can write:
$$\tilde \K_{s,s+\varepsilon}^{s+\varepsilon,\y}=\varepsilon^{-1}\T_\varepsilon \widehat {\tilde \K}_1^{s,s+\varepsilon,s+\varepsilon,\y} \T_\varepsilon ,\ \K_{s,s+\varepsilon}^{s,\x}=\varepsilon^{-1} \T_\varepsilon \hat \K_1^{s,s+\varepsilon,s,\x} \T_\varepsilon , $$
where $\widehat {\tilde \K}_1^{s,s+\varepsilon,s+\varepsilon,\y},\hat \K_1^{s,s+\varepsilon,s,\x} $ are uniformly elliptic and bounded matrices of $\R^{nd}\otimes \R^{nd} $.

Now,  the covariance matrices explicitly write
\begin{eqnarray*}
\tilde \K_{s,s+\varepsilon}^{s+\varepsilon,\y}&=&\int_{s}^{s+\varepsilon} du \tilde \gR^{s+\varepsilon,\y}(s+\varepsilon,u)B a(u,\btheta_{u,s+\varepsilon}(\y))B^* \tilde \gR^{s+\varepsilon,\y}(s+\varepsilon,u)^* ,\\
 \K_{s,s+\varepsilon}^{s,\x}&=&\int_s^{s+\varepsilon} du \gR^{s,\x}(s+\varepsilon,u)B a(u,\btheta_{u,s}(\x))B^* \gR^{s,\x}(s+\varepsilon,u)^* ,
 \end{eqnarray*}
 where $\tilde \gR^{s+\varepsilon,\y},\ \gR^{s,\x} $ respectively denote the resolvents associated to the linear parts of equations \eqref{eq:F:linear:eq:upper} and \eqref{BAR_X}.
Thus, our standing smoothness assumptions in \A{A} (i.e. $a$, $(\nabla_{i-1}\gF_i)_{i\in \leftB 2 ,n\rightB} $ are supposed to be uniformly $\eta $-H\"older continuous) and the bi-Lipschitz property of the flow give:
\begin{eqnarray*}
\left |\langle (\tilde \K_{s,s+\varepsilon}^{s+\varepsilon,\y}-\K_{s,s+\varepsilon}^{s,\x}) (\tilde \btheta_{s+\varepsilon,s}^{s+\varepsilon,\y}(\x)-\y),\tilde \btheta_{s+\varepsilon,s}^{s+\varepsilon,\y}(\x)-\y \rangle\right|\\
=\left| \langle (\widehat{\tilde \K}_1^{s,s+\varepsilon,s+\varepsilon,\y}-\hat \K_1^{s,s+\varepsilon,s,\x})(\varepsilon^{-1/2}\T_\varepsilon (\tilde \btheta_{s+\varepsilon,s}^{s+\varepsilon,\y}(\x)-\y)),\varepsilon^{-1/2}\T_\varepsilon (\tilde \btheta_{s+\varepsilon,s}^{s+\varepsilon,\y}(\x)-\y) \rangle \right| \\
\le C| \btheta_{s+\varepsilon,s}(\x)-\y |^{\eta} \varepsilon^{-1}|\T_\varepsilon (\tilde \btheta_{s+\varepsilon,s}^{s+\varepsilon,\y}(\x)-\y)|^2.
\end{eqnarray*}
Because of the non degeneracy of $a $, the inverse matrices $\left(\widehat {\tilde \K}_1^{s,s+\varepsilon,s+\varepsilon,\y}\right)^{-1}, (\hat \K_1^{s,s+\varepsilon,s,\x})^{-1} $ have the same H\"older regularity. Indeed, up to a change of coordinates one can assume that one of the two matrices is diagonal at the point considered and that the other has dominant diagonal if  $ | \btheta_{s+\varepsilon,s}(\x)-\y |$ is small enough (depending on the ellipticity bound and the dimension). This reduces to the scalar case. Hence,
\begin{eqnarray*}
|Q_\varepsilon|:=\left |\langle \bigl ( (\tilde \K_{s,s+\varepsilon}^{s+\varepsilon,\y})^{-1}-(\K_{s,s+\varepsilon}^{s,\x})^{-1}\bigr) (\tilde \btheta_{s+\varepsilon,s}^{s+\varepsilon,\y}(\x)-\y),\tilde \btheta_{s+\varepsilon,s}^{s+\varepsilon,\y}(\x)-\y \rangle\right|\\
=\left| \langle \bigl(   \left(\widehat{\tilde \K}_1^{s,s+\varepsilon,s+\varepsilon,\y}\right)^{-1}- \bigl(\hat \K_1^{s,s+\varepsilon,s,\x}\bigl)^{-1} \bigl)(\varepsilon^{1/2}\T_\varepsilon^{-1} (\tilde \btheta_{s+\varepsilon,s}^{s+\varepsilon,\y}(\x)-\y)),\varepsilon^{1/2}\T_\varepsilon^{-1} (\tilde \btheta_{s+\varepsilon,s}^{s+\varepsilon,\y}(\x)-\y) \rangle \right| \\
\le C| \btheta_{s+\varepsilon,s}(\x)-\y |^{\eta} \varepsilon |\T_\varepsilon^{-1} (\tilde \btheta_{s+\varepsilon,s}^{s+\varepsilon,\y}(\x)-\y)|^2\le C| \btheta_{s+\varepsilon,s}(\x)-\y |^{\eta} \varepsilon |\T_\varepsilon^{-1} ( \btheta_{s+\varepsilon,s} (\x)-\y)|^2,
\end{eqnarray*}
for $C:=C($\A{A-$\eta$}$)$ using Lemma \ref{eq:F:240409:2}  and the bi-Lipschitz property of the flow $\btheta $ for the last inequality.
From equations \eqref{CTR_PSI_DER}, \eqref{DEF_BTHETA_2} and \eqref{ND_COND}, we eventually get:
 \begin{eqnarray}
|\bTheta_2^\varepsilon(s,\x)|\le C \varepsilon^{\eta/2}\bint{\R^{nd}}^{}\frac{d\y}{\varepsilon^{n^2d/2}}(\varepsilon^{1/2}|\T_\varepsilon^{-1} (\btheta_{s+\varepsilon,s}(\x)-\y) |)^{\eta} \exp(-C\varepsilon|\T_\varepsilon^{-1}(\btheta_{s+\varepsilon,s}(\x)-\y)|^2) \le C\varepsilon^{\eta/2},
\label{CTR_BTH_2}
\end{eqnarray}
for $C:=C($\A{A-$\eta$}$)$. Arguments similar to those employed for  $\bTheta_2^\varepsilon(s,\x) $ can be used to prove $\bTheta_1^\varepsilon(s,\x)\underset{\varepsilon \rightarrow 0}{\longrightarrow} 0  $. The proof then follows from \eqref{CTR_TH_4},\eqref{CTR_R3},\eqref{CTR_BTH_2} recalling the original decomposition \eqref{decomp_bTHETA}. \qed

%%%% Regler le Pb de la notation des constantes.- >On va le faire juste dans la preuve du Lemme, i.e. dans l'appendice.
\bibliographystyle{alpha}
\bibliography{bibli}

\begin{thebibliography}{EPRB99}

\bibitem[BP09]{bass:perk:09}
R.F. Bass and E.A. Perkins.
\newblock A new technique for proving uniqueness for martingale problems.
\newblock {\em From Probability to Geometry (I): Volume in Honor of the 60th
  Birthday of Jean-Michel Bismut}, pages 47--53, 2009.

\bibitem[BPV01]{baru:poli:vesp:01}
E.~Barucci, S.~Polidoro, and V.~Vespri.
\newblock Some results on partial differential equations and asian options.
\newblock {\em Math. Models Methods Appl. Sci}, 3:475--497, 2001.

\bibitem[DM10]{dela:meno:10}
F.~Delarue and S.~Menozzi.
\newblock Density estimates for a random noise propagating through a chain of
  differential equations.
\newblock {\em Journal of Functional Analysis}, 259--6:1577--1630, 2010.

\bibitem[EPRB99]{eckm:99}
J.-P. Eckmann, C.-A. Pillet, and L.~Rey-Bellet.
\newblock Non-equilibrium statistical mechanics of anharmonic chains coupled to
  two heat baths at different temperatures.
\newblock {\em Comm. Math. Phys.}, 201--3:657--697, 1999.

\bibitem[Fri64]{frie:64}
A.~Friedman.
\newblock {\em Partial differential equations of parabolic type}.
\newblock Prentice-Hall, 1964.

\bibitem[HN04]{hera:nier:04}
F.~H\'erau and F.~Nier.
\newblock Isotropic hypoellipticity and trend to equilibrium for the
  fokker-planck equation with a high-degree potential.
\newblock {\em Arch. Ration. Mech. Anal.}, 171--2:151--218, 2004.

\bibitem[H{\"o}r67]{horm:67}
L.~H{\"o}rmander.
\newblock Hypoelliptic second order differential operators.
\newblock {\em Acta. Math.}, 119:147--171, 1967.

\bibitem[IL90]{ishi:lion:90}
H.~Ishii and P.-L. Lions.
\newblock Viscosity solutions of fully nonlinear second-order elliptic partial
  differential equations.
\newblock {\em J. Diff. Equations}, 83:26--78, 1990.

\bibitem[KM00]{kona:mamm:00}
V.~Konakov and E.~Mammen.
\newblock Local limit theorems for transition densities of {M}arkov chains
  converging to diffusions.
\newblock {\em Prob. Th. Rel. Fields}, 117:551--587, 2000.

\bibitem[KM02]{kona:mamm:02}
V.~Konakov and E.~Mammen.
\newblock Edgeworth type expansions for euler schemes for stochastic
  differential equations.
\newblock {\em Monte Carlo Methods Appl.}, 8--3:271--285, 2002.

\bibitem[KM10]{kona:meno:10}
V.~Konakov and S.~Menozzi.
\newblock Weak error for stable driven stochastic differential equations:
  Expansion of the densities.
\newblock {\em To appear in Journal of Theoretical Probability}, 2010.

\bibitem[KMM10]{kona:meno:molc:10}
V.~Konakov, S.~Menozzi, and S.~Molchanov.
\newblock Explicit parametrix and local limit theorems for some degenerate
  diffusion processes.
\newblock {\em Annales de l'Institut Henri Poincar\'e, S\'erie B},
  46--4:908--923, 2010.

\bibitem[MS67]{mcke:sing:67}
H.~P. McKean and I.~M. Singer.
\newblock Curvature and the eigenvalues of the {L}aplacian.
\newblock {\em J. Differential Geometry}, 1:43--69, 1967.

\bibitem[Nor86]{norr:86}
J.~R. Norris.
\newblock Simplified {M}alliavin {C}alculus.
\newblock {\em S\'eminaire de Probabilit\'es}, XX:101--130, 1986.

\bibitem[RBT00]{reyb:thom:00}
L.~Rey-Bellet and L.~Thomas.
\newblock Asymptotic behavior of thermal nonequilibrium steady states for a
  driven chain of anharmonic oscillators.
\newblock {\em Comm. Math. Phys.}, 215--1:1--24, 2000.

\bibitem[Soi94]{soiz:94}
C.~Soize.
\newblock {\em The Fokker-Planck equation for stochastic dynamical systems and
  its explicit steady state solutions}.
\newblock Series on Advances in Mathematics for Applied Sciences, 17. World
  Scientific Publishing Co., Inc., River Edge, NJ, 1994.

\bibitem[SV79]{stro:vara:79}
D.W. Stroock and S.R.S. Varadhan.
\newblock {\em Multidimensional diffusion processes}.
\newblock Springer-Verlag Berlin Heidelberg New-York, 1979.

\bibitem[Tal02]{tala:02}
D.~Talay.
\newblock Stochastic {H}amiltonian dissipative systems: exponential convergence
  to the invariant measure, and discretization by the implicit {E}uler scheme.
\newblock {\em Markov Processes and Related Fields}, 8--2:163--198, 2002.

\end{thebibliography}

\end{document}